\documentclass[12pt]{article}
\usepackage{amsmath, amsthm}
\usepackage{amscd}
\usepackage{amssymb}
\usepackage{array}

\usepackage{color}

\usepackage[normalem]{ulem}

\usepackage{cancel}

\newtheorem{thm}{Theorem}[section]
\newtheorem{theorem}[thm]{Theorem}

\newtheorem{proposition}[thm]{Proposition}
\newtheorem{definition-proposition}[thm]{Definition-Proposition}

\theoremstyle{definition}
\newtheorem{definition}[thm]{Definition}
\newtheorem{example}[thm]{Example}

\newtheorem{observation}[thm]{Observation}

\newtheorem{remark}[thm]{Remark}

\begin{document}

\newcommand{\id}{\relax{\rm 1\kern-.28em 1}}
\newcommand{\R}{\mathbb{R}}
\newcommand{\C}{\mathbb{C}}
\newcommand{\Z}{\mathbb{Z}}
\newcommand{\bp}{\bf p}
\newcommand{\g}{\mathfrak{G}}
\newcommand{\fp}{\mathfrak p}
\newcommand{\fn}{\mathfrak n}
\newcommand{\fh}{\mathfrak h}
\newcommand{\fs}{\mathfrak s}
\newcommand{\fc}{\mathfrak c}
\newcommand{\fk}{\mathfrak k}
\newcommand{\fl}{\mathfrak l}
\newcommand{\e}{\epsilon}

\newcommand{\cU}{\mathcal{U}}
\newcommand{\cY}{\mathcal{Y}}
\newcommand{\cA}{\mathcal{A}}
\newcommand{\cB}{\mathcal{B}}
\newcommand{\cT}{\mathcal{T}}
\newcommand{\cI}{\mathcal{I}}
\newcommand{\cO}{\mathcal{O}}
\newcommand{\cG}{\mathcal{G}}
\newcommand{\cJ}{\mathcal{J}}
\newcommand{\cF}{\mathcal{F}}
\newcommand{\cK}{\mathcal{K}}
\newcommand{\dd}{\mathcal{D}}
\newcommand{\cE}{\mathcal{E}}
\newcommand{\cH}{\mathcal{H}}
\newcommand{\cM}{\mathcal{M}}

\newcommand{\rGL}{\mathrm{GL}}
\newcommand{\rSU}{\mathrm{SU}}
\newcommand{\rSL}{\mathrm{SL}}
\newcommand{\rSO}{\mathrm{SO}}
\newcommand{\rOSp}{\mathrm{OSp}}
\newcommand{\rsl}{\mathrm{sl}}
\newcommand{\End}{\mathrm{End}}
\newcommand{\Hom}{\mathrm{Hom}}
\newcommand{\diag}{\mathrm{diag}}
\newcommand{\rspan}{\mathrm{span}}
\newcommand{\rank}{\mathrm{rank}}
\newcommand{\Gr}{\mathrm{Gr}}
\newcommand{\ber}{\mathrm{Ber}}

\newcommand{\fsl}{\mathfrak{sl}}
\newcommand{\fg}{\mathfrak{g}}
\newcommand{\ff}{\mathfrak{f}}
\newcommand{\fgl}{\mathfrak{gl}}
\newcommand{\fosp}{\mathfrak{osp}}
\newcommand{\bm}{\mathbf{m}}
\newcommand{\fso}{\mathfrak{so}}
\newcommand{\fsu}{\mathfrak{su}}

\newcommand{\str}{\mathrm{str}}
\newcommand{\Sym}{\mathrm{Sym}}
\newcommand{\tr}{\mathrm{tr}}
\newcommand{\defi}{\mathrm{def}}
\newcommand{\Ber}{\mathrm{Ber}}
\newcommand{\spec}{\mathrm{Spec}}
\newcommand{\sschemes}{\mathrm{(sschemes)}}
\newcommand{\sspaces}{\mathrm{(sspaces)}}
\newcommand{\smflds}{\mathrm{ {(smflds)} }}
\newcommand{\mflds}{\mathrm{ {(mflds)} }}
\newcommand{\sets}{\mathrm{(sets)}}
\newcommand{\mfld}{\mathrm{ {(mfld)} }}
\newcommand{\odd}{\mathrm{odd}}
\newcommand{\alg}{\mathrm{alg}}
\newcommand{\salg}{\mathrm{(salg)}}
\newcommand{\diff}{\mathrm{diff}}
\newcommand{\wa}{\mathrm{(wa)}}
\newcommand{\swa}{\mathrm{(swa)}}
\newcommand{\Amflds}{\mathrm{(\cA_0mflds)}}
\newcommand{\spts}{\mathrm{(spts)}}
\newcommand{\ad}{\mathrm{ad}}
\newcommand{\Ad}{\mathrm{Ad}}
\newcommand{\sloc}{\mathrm{(sloc)}}
\newcommand{\Lie}{\mathrm{Lie}}
\newcommand{\Proj}{\mathrm{Proj}}
\newcommand{\uspec}{\underline{\mathrm{Spec}}}
\newcommand{\uproj}{\mathrm{\underline{Proj}}}
\newcommand{\rOsp}{\mathrm{Osp}}
\newcommand{\rSp}{\mathrm{Sp}}
\newcommand{\rodet}{\mathrm{odet}}
\newcommand{\amflds}{\mathrm{($A$-manifolds)}}
\newcommand{\azmflds}{\mathrm{({A_0}-manifolds)}}

\newcommand{\sym}{\cong}

\newcommand{\al}{\alpha}
\newcommand{\be}{\beta}
\newcommand{\lam}{\lambda}
\newcommand{\de}{\delta}
\newcommand{\ga}{\gamma}
\newcommand{\ep}{\epsilon}
\newcommand{\lra}{\longrightarrow}
\newcommand{\ra}{\rightarrow}
\newcommand{\La}{\Lambda}

\newcommand{\nil}[1]{\accentset{\circ}{#1}} 


\bigskip

\noindent
\centerline{\Large{ \bf Quotients in supergeometry}}

\medskip


\bigskip

\centerline{ L. Balduzzi$^\natural$,
C. Carmeli$^\natural$, R. Fioresi$^\flat$}

\bigskip

\centerline{\it $^\natural$ Dipartimento di Fisica, Universit\`a di Genova 
and INFN, sezione di Genova}
\centerline{\it Via Dodecaneso, 33 16146 Genova, Italy}
\centerline{
{\footnotesize e-mail: 
luigi.balduzzi@ge.infn.it, claudio.carmeli@ge.infn.it}}
\bigskip

\centerline{\it $^\flat$ Dipartimento di Matematica,
              Universit\`a di Bologna }
\centerline{\it Piazza di Porta San Donato, 5
40127 Bologna, Italy}
\centerline{{\footnotesize e-mail: fioresi@dm.unibo.it}}

\begin{abstract} The purpose of this paper is to present the notion of
quotient of supergroups in different categories using the
unified treatment of the functor of points and to
examine some physically interesting examples.
\end{abstract}

\section{Introduction}

The study of supergeometry was prompted by important physical
questions linked to the symmetries of physical systems, which
take into account the intrinsecally different nature of the two
fundamental types of particles: bosons and fermions.

\medskip

While the bosons obey the Bose-Einstein statistics,
the fermions are described by the Fermi one. These two types of
particles have a fundamentally different behaviour:
the bosons are described by \textit{commuting} functions,
while the fermions by \textit{anticommuting} ones. Since these
particles do transform into each other, it is necessary
to consider symmetries which allow to mix
these two types of functions.

\medskip

From a purely mathematical point of view, we can view supergeometry
as  $\Z_2$-graded geometry,
where every ordinary geometric concept, as for example
manifolds, varieties, vector fields and so on,
has an $\Z_2$-graded corresponding
one. It is however important to stress
that a supermanifold is not to be understood as an ordinary
manifold with an associated $\Z_2$-graded vector bundle,
since in supergeometry we allow transformations which
mix the even and the odd coordinates, as we shall see in Section
\ref{prelim}.

\medskip

Our treatment is organized as follows.

\medskip

In Section \ref{prelim} we quickly review some general facts
on supergeometry including the functor of points
approach to the study of superspaces.

\medskip

In Section \ref{homospaces} we define what an action of
a supergroup on a superspace is and the concept of
homogeneous superspace.

\medskip

In Section \ref{functorofpoints} we
define the functor of points and the functor of $A$-points for
homogeneous spaces. We also examine in detail the example
of the superflag and  its big cell, together with its
physical interpretation as superconformal and super Minkowski
spaces.


\medskip

We want to especially thank prof. V. S. Varadajan
for his constant encouragement and his generosity in sharing his time and
his ideas with us at all times and also while preparing this paper.

\section{Preliminaries} \label{prelim}

Let $k$ be the ground field, $char(k) \neq 2,3$.

\medskip

For the basic definitions of superalgebra, supervector
space and similar, refer to \cite{vsv} ch. 4 and \cite{ma1}
ch. 3

\medskip

\begin{definition}

A \textit{superspace} $S=(|S|,\cO_S)$ consists of a
topological space $|S|$ together with a sheaf of
commuting superalgebras
$\cO_S$, with the property that the
stalk $\cO_{S,x}$ is a local superalgebra for all $x \in |S|$.
A morphism of superspaces $\varphi:S \longrightarrow T$ is a
continuous map $|\varphi|: |S| \longrightarrow |T|$ together with a sheaf map
$\varphi^*:\cO_T \longrightarrow \varphi_*\cO_S$
so that $\varphi^*_x(\mathfrak{m}_{|\varphi|(x)}) \subset
\mathfrak{m}_x$ where
$\mathfrak{m}_x$ is the maximal ideal in
$\cO_{S,x}$ and $\varphi^*_x$ is the
stalk map.
by $\varphi: S \longrightarrow T$.

\end{definition}

We shall denote with $\sspaces$ the category of superspaces.

\medskip
Let's see some key examples of superspaces.

\begin{example}

1. $\R^{p|q}$. On the topological space $\R^p$ we define the sheaf
of commutative $\R$-superalgebras:
\[
V \mapsto \cO_{\R^{p|q}}(V):=C_{\R^p}^{\infty}(V)[\theta_1, \ldots, \theta_q],
\]
where $C_{\R^p}^{\infty}(V)[\theta^1, \ldots, \theta^q]=
C_{\R^p}^{\infty}(V) \otimes \wedge(\theta_1, \ldots, \theta_q)$ and
the $\theta_j$ have to be thought as odd (anti-commuting) indeterminates.

One can readily check that $\R^{p|q}:=(\R^p,\cO_{\R^{p|q}})$ is a
superspace. Notice that the morphisms of superspaces are
allowed to mix even and odd coordinates. For example we can define
the morphism $\phi:\R^{1|2} \lra \R^{1|2}$ on global
section by: $\phi(x)=x+\theta_1\theta_2$, $\phi(\theta_1)=\theta_1$,
$\phi(\theta_2)=\theta_2$. This tells that $\R^{p|q}$ cannot be
simply viewed as $\R$ together with an exterior bundle.

\medskip\noindent
2. $\R^{p|q}_h$, $\C^{p|q}_h$.
Similarly define for $V \subset \R^p$ open, the sheaf of superalgebras:
\[
V \mapsto \cH_{\R^{p|q}}(V):=\cH_{\R^p}(V)[\theta_1, \ldots, \theta_q]:=
\cH_{\R^p}(V) \otimes \wedge(\theta_1, \ldots, \theta_q).
\]
where $\cH_{\R^p}$ denotes the sheaf of real analytic functions on $V$.
Again one can check that $\R^{p|q}_h=(\R^p,\cH_{\R^{p|q}})$
is a superspace. The definition of the
superspace $\C^{p|q}_h=(\C^p,\cH_{\C^{p|q}})$ goes along the
same lines.

\medskip\noindent
3. $\uspec A$. \label{spec}
Let $A$ be a commutative superalgebra. Since $A_0$ is an algebra, we can
consider the topological space $$\spec(A_0)=\{\hbox{prime ideals }
\bp\subset A_0\}.$$
The closed sets are $V(S)=\{\fp \in \spec(A_0) \quad | \quad
\fp \supset S\}$. Classically we can define the
structural sheaf $\cO_{A_0}$ on $\spec(A_0)$
by giving on an open cover of $\spec(A_0)$ by $U_i=\spec(A_0[f_i^{-1}])$
the sheaves $\cO_{A_0}|_{U_i}(U_i):=A_0[f_i]$. The stalk
of the structural sheaf at the prime $\bp\in \spec(A_0)$ is the
localization of $A_0$ at $\bp$. We can replicate this construction in
the super setting.
As for any superalgebra, $A$ is a module over $A_0$, and we have
indeed a sheaf $\tilde A$ of $\cO_{A_0}$-modules over $\spec A_0$
with stalk $A_{\bp}$, the localization of the $A_0$-module $A$ over
each prime $\bp \in \spec(A_0)$.
$\uspec A=_{\defi}(\spec A_0,\tilde A)$ is a superspace. As before
$\uspec A$ is covered by open subsuperspaces $U=\uspec A[f^{-1}]$,
$f \in A_0$.
(For more details concerning  the
construction of the sheaf $\tilde M$ for a generic $A_0$ module $M$,
see Ref. \cite{ha} II \S 5 and \cite{eh} Ch. 1).

\end{example}

\begin{definition}
We say that a superspace
$M$ is a \textit{supermanifold} (resp. \textit{real or complex analytic
supermanifold})
if $M$ is locally isomorphic to $\R^{p|q}$ (resp. $\R^{p|q}_h$ or
$\C^{p|q}_h$).
We also say that a superspace
$M$ is a \textit{superscheme} if it is
locally isomorphic to the spectrum of some superalgebra
(of course the superalgebras may be different at different points).
\end{definition}

\begin{definition}
Given a superspace $G$, if we have three
morphisms:
$$
m:G \times G \lra G, \qquad i:G \lra G, \qquad 1:\{\bullet\} \lra G
$$
satisfying the usual commutative diagrams for
multiplication, inverse and identity in
an abstract group, we say that $G$ is a \textit{supergroup}. If
furtherly $G$ is a supermanifold, (resp. complex or real analytic),
we say $G$ is a \textit{Lie} (resp. \textit{complex or real
analytic}) supergroup. If $G$ is a superscheme, we say that $G$ is a
\textit{supergroup scheme}.
\end{definition}

The concept of functor of points allows us to recover some
of the geometric intuition.

\begin{definition} \label{fopts}
We define the \textit{functor of points} $h_X$ of the superspace $X$
as the representable functor
$$
h_X:\sspaces \lra \sets, \qquad T \mapsto h_X(T)=\Hom(T,X).
$$
In the same way, by the appropriate changes
in the categories, we can define the functor of points of a
supermanifold or a superscheme.
Clearly if the superspace $G$ is a supergroup, the functor is group-valued
(and vice-versa).
\end{definition}

The functor of points approach is so powerful because of Yoneda's
Lemma, that we state in a special form of interest to us:


\begin{theorem} Yoneda's Lemma. We have a bijection between the set
of morphisms of supermanifolds (supervarieties) $X \lra Y$
and the set of natural transformations $h_X \lra h_Y$.
\end{theorem}

\begin{observation} \label{sheafification}
By its very definition the functor of points $h_S$ of a superspace
$S$ has the \textit{presheaf property}, that is, when restricted to
the open subsets of a superspace it is a presheaf of sets (recall that
a presheaf is just a functor from the category of open sets
of a topological space, where the morphisms are given by inclusions).
However $h_S$
has also the \textit{sheaf} property; in other words if $\{T_i\}$ is
a covering of the superspace $T$ and we have a family $\al_i \in
h_S(T_i)$, such that $\al_i|_{T_i \cap T_j}=\al_j|_{T_i \cap T_j}$,
then there exists a unique $\al \in h_S(T)$ such that
$\al|_{T_i}=\al_i$\footnote{
As customary we denote $\al|_{T_i}$ as the image of
$\al \in h_S(T)$ under the map $h_S(\phi_i)$, where
$\phi_i: T_i \hookrightarrow T$.}. We leave this verification as an exercise to
the reader.

\medskip

Any functor $F:\sspaces \lra \sets$ is a presheaf and
as, for any presheaf, we can always build its \textit{sheafification}
$\widetilde{F}: \sspaces \lra \sets$, which has the
following properties:
\\
1. $\tilde F$ is a sheaf.
\\
2. There is a canonically defined presheaf
morphism $\psi:F \lra \tilde F$.
\\
3. Any presheaf morphism $\phi:F \lra G$, with $G$ sheaf, factors
via $\psi$, i. e. 
$\phi:F \stackrel{\psi} \lra
\tilde F \lra G$.
\\
Moreover $\tilde F$ is
locally is isomorphic to $F$. For more
details on this construction we refer the reader to
\cite{dg} and \cite{eh}.

\end{observation}
\medskip

Next, we want to introduce the concept of $\cA_0$-manifold
and the functor of the $A$-points of a supermanifold $X$.
This is substantially different from the functor of points $h_X$ we
have already described; in fact we can define it only in the
differential and holomorphic categories. 
We are going to see that it characterizes the
supermanifold and in many computational problems it allows to
simplify significantly the notation. For a complete treatment see
\cite{bcf}.

\medskip

Let our ground field $k$ be $\R$ or $\C$.

\begin{definition}
We call the commutative algebra $A$ 
a \textit{Weil algebra} if it is local, finite dimensional and $A=k \oplus J$,
with the nilpotent maximal ideal $J$.
 We denote with $\wa$ the category of Weyl algebras
(sometimes called \textit{local algebras}) and with $\swa$ the
category of Weyl superalgebras, defined in a similar way.


Let $A_0$ be a local algebra (the
index $0$ reminds us it has no odd elements). A manifold $M$ is
called an $A_0$\textit{-manifold} if there is an $A_0$-module
$L$ and an open cover $\{U_i\}$ of $M$, such that 
$h_i: U_i \lra U'_i \subseteq L$ are
diffeomorphisms (of $C^\infty$ manifolds) and $d(h_i \cdot
h_j^{-1})$ are isomorphisms of $A_0$-modules. The set of all
$A_0$-manifolds for all $A_0 \in \wa$ forms the objects of the
category of $\cA_0$-manifolds that we denote with ${\Amflds}$. A
morphism of two $\cA_0$-manifolds $M$ and $N$, $M$ being an
$A_0$-manifold, $N$ a $B_0$-manifold, consists of a pair $(f,\phi)$,
where $f:M \lra N$ is $C^\infty$ morphism and $\phi:A_0 \lra B_0$ an
algebra morphism such that $df(ax)=\phi(a) df(x)$.
\end{definition}

We are ready to define the functor of the $A$-points
of a supermanifold, through a definition-proposition (more
details can be found in \cite{bcf}).



\begin{definition-proposition}
Let $M$ be a supermanifold. We define the {\sl set of $A$-points} of
$M$
$$
M_A := \coprod_{x\in |M|}\Hom_\salg(\cO_{M,x}, A)
$$
It has a natural structure of $A_0$-manifold. We define the {\sl
local functor of points} of $M$ the functorial assignment
\begin{eqnarray*}
M_{(.)}:\swa \lra \Amflds, \qquad A \mapsto M_A .
\end{eqnarray*}
\end{definition-proposition}

For more details see \cite{bcf}.
\medskip


When $M$ is smooth, we can write the functor $M_A$ in a
much simpler way (see \cite{bcf}).

\begin{proposition}
Let $M$ be a smooth supermanifold, then:
$$M_A \cong \Hom(\cO_M(M),A).$$
\end{proposition}


As it happens for the functor of points $h_X$, also in this case
we can give an analogue of Yoneda's lemma. 
This means that the the functor $\cY$, $\cY(M)=_{\defi}M_{(.)}$
is  a fully faithful embedding.
As for the usual functor of points,
$\cY$ is not an equivalence of categories.
In other words, not all the functors $h:\, \swa \to \Amflds$
arise as  the functors of $A$-points of a  super manifold.
If this is the case, in analogy with the functor of points notation, we say
the functor is {\textit{ representable}}.
In this frameworks it is possible to prove the
following representability criterion, that we state for both
the functor of $A$-points and the functor of points discussed in
\ref{fopts}.

\begin{proposition} 
1. Let $F: \smflds \lra \sets$ be a functor with the sheaf property.
Suppose that $F$ admits
a cover by open subfunctors, i. e. there exist
representable subfunctors of $F$, $U_i:\smflds \lra \sets$, such that for
any supermanifold $M$ and any natural transformation $f:h_M
\lra F$, $f^{-1}(U_i)=h_{V_i}$ and the $V_i$ are open and cover $M$. 
Then $F$ is
representable, i. e. it is the functor of points of
a supermanifold.\\ 
2. Let $h:\, \swa \to \Amflds$ be a functor. Denote by $p_A:A \lra
\R$ the canonical projection of an algebra $A \in \swa$ into $A/J \cong \R$.
Suppose that an open cover
$\{ \widetilde{U}_\alpha\}$ of $h(\R^{0|0})$
is given such that the functors
\begin{align*}
h_\alpha:\,  \swa &\to \Amflds & 
A & \mapsto (h_{p_A})^{-1}(\widetilde{U}_\alpha) 
\end{align*}
are representable by $\R^{n|m}$, for fixed $n$ and $m$. Then
$h$ is representable, i. e. it is the functor of the $A$-points
of a supermanifold.
\end{proposition}
\begin{proof}
For (1) see \cite{flv}, for (2) see \cite{bcf}.
\end{proof}

As we shall see in the next sections, this is an important
result that allows us to define properly the quotients of
supergroups and their functor of points.

\medskip

\section{Actions of supergroups on superspaces} \label{homospaces}

Let $k$ be the ground field, $char(k) \neq 2,3$

\begin{definition} \label{action}
Let $G$ be a supergroup. We say that $G$ \textit{acts} on the
superspace $M$ if there exists a morphism $\phi:G \times M \lra M$
denoted as $(g,x)  \mapsto  g \cdot x$ for $g \in G(T)$ and $x \in
M(T)$,
such that for all superspaces $T$: \\
1. $1 \cdot x=x$, $\forall x \in M(T)$ \\
2. $(g_1g_2) \cdot x=g_1 \cdot (g_2 \cdot x)$, $\forall x \in M(T)$,
$\forall g_1, g_2 \in G(T)$.

\medskip

We say that $G$ acts \textit{transitively} on $M$, or that $M$ is an
\textit{homogeneous space} if there is $x_0 \in |M|$ such that the
morphism $\phi_{x_0}:G \lra M$, $\phi_{x_0}(g)=g \cdot x$ is onto,
i.e. the sheafification ${\widetilde {Im(\phi_{x_0})}}$
of the image presheaf coincides with $M$ (see \ref{sheafification}).
\end{definition}

\medskip

One can give in an obvious way this same definition in the categories
of supermanifolds and superschemes.

\medskip

When $M$ is a supermanifold, our definition
of homogeneous space
is equivalent to the one appearing in \cite{bcc} as the
next proposition shows.

\begin{theorem} \label{rep}
$\widetilde{Im \phi_{x_0}} = M$ if and only if $\phi_{x_0}$ is a
surjective submersion.
\end{theorem}

\begin{proof}
For brevity let $\phi = \phi_{x_0}$. Let us suppose that $\phi$ is a
surjective submersion. Let $m \in |M|$ and $g \in
|\phi|^{-1}(m)$ ($|\phi|$ is surjective, so it exists).
Since $\phi$ is a submersion there exists $V \subseteq |G|$
with coordinates $X_1,\ldots,X_{p+q}$ ($dim G = p|q$) and $W
\subseteq |M|$ with coordinates $Y_1,\ldots,Y_{m+n}$ ($dim M =
m|n$) such that
\[
    \phi^*(Y_i) = X_i
\]
Let $t \in U \subseteq |T|$ and $\alpha \colon U \to M$ such
that $m = |\alpha|(t)$. We can suppose $|\alpha|(U)
\subseteq W$. If $\alpha^*(Y_i) = f_i \in \cO_T(U)$, $\beta
\colon U \to V$ defined by
\[
    \beta^*(X_i) =
    \begin{cases}
        f_i &\text{if $i \leq m+n$} \\
        0 &\text{otherwise}
    \end{cases}
\]
satisfies $\phi \circ \beta = \alpha$. Then $[\alpha] \in (Im
\phi)_t$, hence $(Im\phi)_t=M_t$ and this gives one implication.

Vice-versa let us suppose that $\widetilde{Im \phi} = M$. Taking $T
= \R^{0|0}$ we have that $|\phi|$ must be surjective. Let's  now
assume $T = M$ and $m \in |M|$. There exists $U \ni m$ and $\psi
\colon U \to G$ such that $\phi \circ \psi = \id_U$. Then $\phi$
must be a submersion at $|\psi|(m)$ and this is true
everywhere, since $\phi$ has constant rank. Indeed for all $g \in
|G|$,
$$
{(d\phi)}_g \circ {(dl^G_g)}_1 = {(dl^M_g)}_{x_0} \circ {(d\phi)}_1
$$
where the isomorphisms $l^G_g$ and $l^M_g$ are the left actions of
$g$ on $G$ and $M$ respectively. 
\end{proof}

%
%
%
%

\begin{definition} Let's the notation be as above.
The functor:
$$
S_{x_0}(T)=\{ g \in G(T) \quad | \quad g \cdot x_0=x_0 \}, \quad T
\in \sspaces 
$$
is called the \textit{stabilizer} of $x_0 \in |M|$.
\end{definition}

We have given this definition in general, however we are especially
interested in two cases:
\\
1. $G$ Lie supergroup, $M$ a supermanifold.
\\
2. $G$ complex algebraic supergroup, $M$ complex algebraic variety.
\\
In each case the definitions above need to be suitably modified taking the
superspaces in the appropriate category.

\medskip


\begin{theorem} \label{stabthm}
Let $G$ be a Lie or algebraic affine supergroup acting transitively
on the supermanifold or supervariety $M$, $x_0
\in |M|$. Then
\\ \\
1. $S_{x_0}^{\diff}: \smflds \lra \sets$,
$S_{x_0}^\diff(T)=\{ g \in G(T) \quad | \quad g \cdot x_0=x_0 \}$,
\\ \\
2.  $S_{x_0}^{\alg}: \salg \lra \sets$,
$S_{x_0}^\alg(A)=\{ g \in G(A) \quad | \quad g \cdot x_0=x_0 \}$, \\ \\
are the functor of points respectively of a Lie supergroup and
of an algebraic supergroup. In other words the stabilizer supergroup
functor is representable.
\end{theorem}

\begin{proof}
For the differential category see \cite{dm} and \cite{bcc},
while for the algebraic category, see \cite{fi-smooth}.
\end{proof}



\medskip

There are many examples of actions of supergroups on superspaces,
some of which are especially interesting.
We now are going to see that Theorem \ref{stabthm} gives the representability
for all the classical
supergroups both in the categories of Lie and algebraic
supergroups.

\medskip

Let $k$ be the field $\R$ or $\C$ for the supermanifolds
category and just a generic field, with $char(k) \neq 2,3$
for the superschemes category.

\medskip

1. {\bf $A(n)$ series}.
Let's first consider the algebraic setting.
Let $A \in \salg$.
Define $\rGL_{m|n}(A)$ as the set of all invertible
morphisms $g:A^{m|n}\rightarrow A^{m|n}$.
This is equivalent to ask
that the {\it Berezinian} \cite{be} or {\it superdeterminant}
$$
\Ber(g)=\Ber \begin{pmatrix} p & q \\ r & s \end{pmatrix} =
\det(p-qs^{-1}r)\det(s^{-1})
$$
is invertible in $A$
(where $p$ and $s$ are $m \times m$, $n \times n$ matrices
of even elements in $A$, while $q$ and $r$ are $m \times n$,
$n \times m$ matrices of odd elements in $A$).
A necessary and sufficient
condition for $g \in \rGL_{m|n}(A)$ to be invertible is that
$p$ and $s$ are invertible. The group valued functor
$$
\begin{array}{cccc}
\rGL_{m|n} &:\salg & \longrightarrow & \sets \\ 
&A&\longmapsto &
\rGL_{m|n}(A).
\end{array}
$$
is an affine supergroup called the \textit{general
linear supergroup} and it is represented by the algebra
\begin{eqnarray*}&
k[\rGL_{m|n}]:=k[x_{ij}, y_{\alpha \beta},
\xi_{i\beta},\gamma_{\alpha
j},z,w]/\bigr((w\det(x)-1,z\det(y)-1\bigl),
\\
&i,j=1,\dots m,\;\;
\alpha,\beta=1,\dots n.
\end{eqnarray*}

Consider the morphism

\begin{eqnarray} \label{A}
\rho: \rGL_{m|n} \times k^{1|0} \lra k^{1|0} \qquad
(g,c) \lra \Ber(g)c.
\end{eqnarray}

The stabilizer of the point $1 \in k^{1|0}$
coincides with all the matrices in $\rGL_{m|n}(A)$ with
Berezinian equal to 1, that is $\rSL_{m|n}(A)$ the
special linear supergroup. By the Theorem \ref{stabthm}
we have immediately that $\rSL_{m|n}$ is representable
as an algebraic supergroup.

\medskip
The supermanifold case is very similar. Define the functor (by an
abuse of notation we use the same symbol) $\rGL_{m|n}(T)$ as the
invertible $\cO_T$-module sheaf morphisms $\cO_T^{m|n} \lra
\cO_T^{m|n}$. $\rGL_{m|n}(T)$ can also be identified with the $m|n$
matrices with coefficients in $\cO_T(|T|)$. In fact any morphism of
supermanifold sheaves is determined once we know the morphism on the
global sections $\cO_T^{m|n}(|T|) \lra \cO_T^{m|n}(|T|)$. Again we
can define the Berezinian of a matrix and we can consider a morphism
as in \ref{A}. The stabilizer of the point $1 \in k^{1|0}$ coincides
with all the matrices in $\rGL_{m|n}(T)$ with Berezinian equal to 1,
that is $\rSL_{m|n}(T)$ the special linear Lie supergroup. By the
Theorem \ref{stabthm} we have that $\rSL_{m|n}$ is representable as
a Lie supergroup.

\medskip

2. {\bf $B(m,n)$, $C(n)$, $D(m,n)$ series}.
Consider the morphism (both in the superscheme and
supermanifold categories):
\begin{eqnarray} \label{BCD}
\rho: \rGL_{m|2n} \times \cB \lra \cB \qquad
(g,\psi(\cdot,\cdot)) \lra \psi(g \cdot, g \cdot),
\end{eqnarray}
where $\cB$ is the supervector space of all the
symmetric bilinear forms on $k^{m|2n}$.
We define $\rOsp_{m|2n}$ as
the stabilizer of the point $\Phi$, the standard
bilinear form on $k^{m|2n}$. Again this is an algebraic and Lie supergroup
by Theorem \ref{stabthm}.

\medskip

3. $P(n)$ {\bf series}.
Define the algebraic and Lie supergroup $\pi \rSp_{n|n}$ as we did
for $\rOsp_{m|n}$, by taking antisymmetric bilinear forms
instead of symmetric ones. Consider the action:
$$
\begin{array}{ccc}
\pi \rSp_{n|n} \times k^{1|0}  \lra  k^{1|0} \qquad
(g,c)  \mapsto  \Ber(g)c.
\end{array}
$$
By Theorem \ref{stabthm} we have that $Stab_1$ is an affine
algebraic supergroup, hence it is an algebraic and
Lie supergroup. It is corresponding
to the $P(n)$ series.

\medskip

3. $Q(n)$ {\bf series}. Let $D=k[\eta]/(\eta^2+1)$. This is a
non commutative superalgebra. Define the supergroup functor
$GL_n(D):\salg \lra \sets$, with $GL_n(D)(A)$ the group  of automorphisms
of the left supermodule $A \otimes D$. In \cite{dm} is proven
the existence of a morphism called the \textit{odd determinant}
$$
\rodet: GL_n(D) \lra k^{0|1}.
$$
Reasoning as before define:
$$
GL_n(D) \times k^{0|1} \lra k^{0|1}, \qquad g,c \lra \rodet(g)c.
$$
Then $G=Stab_1$ is an affine algebraic supergroup and for $n \geq 2$
we define $Qg(n)$ as the quotient of $G$ and the diagonal subgroup
$GL_{1|0}$. This is an algebraic and Lie supergroup
and its Lie superalgebra is $Q(n)$.

\section{Homogeneous spaces via their functor of points}
\label{functorofpoints}

We now want to address the following question. Let $G$ be
a supergroup and $H$ a closed subgroup, i. e.
$|H|$ is closed in $|G|$. Consider the functor:
$$
\sspaces \lra \sets, \qquad T \mapsto G(T)/H(T).
$$
Is this functor representable?
In this generality the answer is no,
however we shall describe a representability result in the
categories of supermanifolds and supervarieties.

\begin{theorem}
Let $G$ be a Lie supergroup, $H$ a closed Lie
subgroup. Let $\widetilde{G/H}$ be the
sheafification of the functor:
$$
T \mapsto G(T)/H(T).
$$
Then $\widetilde{G/H}$ is the functor of points of a
supermanifold that we denote with $G/H$. Moreover $G/H$ is unique
supermanifold with underlying topological space $|G|/|H|$
with respect to the following property:
\\
The natural morphism $\pi:G \lra G/H$ is a submersion, moreover
$G$ acts on $G/H$ and we have the commutative diagram:
$$
\begin{array}{ccc}
G \times G & \stackrel{m}  \lra & G  \\
\downarrow & & \downarrow \pi \\
G \times G/H & \lra & G/H
\end{array}
$$
\end{theorem}

\begin{proof} A complete proof of this statement can be
found in \cite{flv}.
\end{proof}

\begin{remark}
In the algebraic setting, Zubkov recently proved in \cite{zu} a similar result
for $G/H$ affine and in the case of $char(k)=0$. In this
setting one has to be more careful in taking the sheafification
and more difficulties are present, since we don't have in
general the local splitting of $G$ as $H \times W$ at the identity.
\end{remark}

We now turn to the formulation of the same problem for the
functor of the $A$-points.

\begin{proposition}
Let $G$ be a Lie supergroup and $H$ be a closed subgroup. The functor
\begin{align*}
\swa &\to \Amflds\\
A &\mapsto G_A/H_A
\end{align*}
is representable.
\end{proposition}

\begin{proof}
It is well known
that there exists an open cover of $G$ by tubular neighborhoods
$U_\alpha \cong W_\alpha \times H$,
where $W_\alpha$ are isomorphic to open sub superdomains in $\R^{p|q}$.
Since the functor of $A$-points is product preserving we have that
\begin{eqnarray*}
(U_\alpha)_A/H_A \cong (W_\alpha)_A
\end{eqnarray*}
and the result follows immediately from the Representability Theorem
\ref{rep}.
\end{proof}

As an example, we shall examine the construction of the
superflag $\cF$ of $2|0$ and $2|1$ spaces in the $4|1$ dimensional
complex super vector space $\C^{4|1}$.
This is important in physics,
since it gives the complexification of the super conformal
space containing as big cell the Minkowski superspace
(for more details on the physical interpretation see
\cite{flv}).


\medskip

Let $\cF$ be the functor:
$\cF: \smflds \lra \sets$, where $\cF(T)$ is the
set of $2|0$ and $2|1$ projective modules $Z_1 \subset Z_2$ inside
$\cO_T^{4|1}:=\cO_T \otimes \C^{4|1}$.
$\cF$ is the functor of points of a supermanifold called the \textit{superflag}
of $2|0$ and $2|1$ planes in $\C^{4|1}$, that we shall
still denote by $\cF$ by an abuse of notation.  Clearly $\cF \subset
\cG_1 \times \cG_2$, where $\cG_1$ and $\cG_2$ are
respectively the supergrassmannians of $2|0$ and $2|1$ planes
in $\C^{4|1}$ (for a direct proof of the non trivial fact that $\cF$, $\cG_1$,
$\cG_2$ are supermanifolds see \cite{ma1}).

\medskip

We are now going to realize $\cF$ as the quotient of $\rSL_{4|1}$
by a suitable parabolic subgroup.

\medskip

The natural action of $G=\rSL_{4|1}$ on $\cO_T^{4|1}$
induces an action on $\cG_1$ and $\cG_2$ and also on $\cF$:
$$
\begin{array}{cccc}
G(T) & \lra & \cF(T) \subset & \cG_1(T) \times \cG_2 \\
g & \mapsto & g \cdot F. &
\end{array}
$$

Let us fix the element
$F_0=\{\cO_T^{2|0} \subset \cO_T^{2|1}\}$ 
in $\cF(T)$. 
Then we can write the action as:
$$
\begin{array}{c}
g \cdot F_0=\left(
\begin{pmatrix}
g_{11} & g_{12} \\
g_{21} & g_{22} \\
g_{31} & g_{32} \\
g_{41} & g_{42} \\
\gamma_{51} & \gamma_{52}
\end{pmatrix},
\begin{pmatrix}
g_{11} & g_{12} & \gamma_{15}\\
g_{21} & g_{22} & \gamma_{25} \\
g_{31} & g_{32} & \gamma_{35} \\
g_{41} & g_{42} & \gamma_{45}\\
\gamma_{51} & \gamma_{52} & g_{55}
\end{pmatrix} \right) \in \cG_1(T) \times \cG_2(T).
\end{array}
$$

The stabilizer subgroup functor at $F_0$ is given as the subgroup
$H(T)$ of $G(T)$
consisting of all matrices in $G(T)$ of the
form:
$$
\begin{pmatrix}
g_{11} & g_{12} & g_{13} & g_{14} & \gamma_{15} \\
g_{21} & g_{22} & g_{23} & g_{24} & \gamma_{25} \\
0 & 0 & g_{33} & g_{34} & 0 \\
0 & 0 & g_{43} & g_{44} & 0 \\
0 & 0 & \gamma_{53} & \gamma_{54} & g_{55} \\
\end{pmatrix}.
$$
$H$ is clearly representable by a group supermanifold 
moreover we have that locally:
$$
T \mapsto G(T)/H(T)=\cF(T).
$$
Hence $\widetilde {G/H}=\cF$
and this is the functor of points of the superflag $\cF=G/H$.

\medskip

We wish now to describe explicitly $G/H$ and its {\sl big cell} $U$
and to prove explicitly that the map $\pi:G \lra G/H$ is
a submersion.


The big cell $U$ in $\cF$ is defined as $\cF \cap v_1 \times v_2$,
where $v_1$ and $v_2$ are the big cells inside $\cG_1$ and $\cG_2$.
By definition
$v_1(T)$ contains all the elements in $\cG_1(T)$ having the determinant
in the upper left corner invertible, while $v_2(T)$ contains all the
elements in $\cG_2(T)$ having the berezinian of rows $1,2,5$ and
columns $1,2,3$ invertible. Hence we can write:
$$
v_1(T)=\begin{pmatrix}
I_2 \\
A \\
\alpha
\end{pmatrix}, \qquad
v_2(T)=\begin{pmatrix}
I_2 & 0 \\
B & \beta \\
0 & 1
\end{pmatrix}  \qquad T \in
\smflds,
$$
where $I_2$ is the identity matrix, $A$ and $B$ are $2 \times 2$ matrices
with even entries and $\al=(\al_1, \al_2)$, $\be^t=(\be_1,\be_2)$
are rows with odd entries.

An element of $v_1(T)$  is inside $v_2(T)$ if and only if
\begin{equation}
A=B+\beta\alpha,\label{twistor}
\end{equation} so we can take as coordinates for a
flag in the big cell $U$ the triplet $(A,\alpha,\beta)$. We see
then  that $U$ is an affine $4|4$ superspace. Equation
(\ref{twistor}) is also known as {\it twistor relation},
in the physics literature.

In these coordinates, $F_0=\left(\begin{pmatrix}
I \\
0 \\
0
\end{pmatrix}, \;
\begin{pmatrix}
I & 0 \\
 0& 0 \\
0 & 1
\end{pmatrix} \right)
$
is described by $(0,0,0)$.

We want to write the map $\pi$ in these coordinates. In a suitable
open subset near the identity of the group we can take an element
$g\in G(T)$ as
$$g=\begin{pmatrix}g_{ij}&\gamma_{i5}\\\gamma_{5j}&g_{55}\end{pmatrix},
\qquad i,j=1,\dots 4.$$ Then, we can write an element $g \cdot \cF
\in \cG_1 \times \cG_2$ as:
\begin{equation}
\begin{pmatrix}
g_{11} & g_{12} \\
g_{21} & g_{22} \\
g_{31} & g_{32} \\
g_{41} & g_{42} \\
\gamma_{51} & \gamma_{52}
\end{pmatrix},
\begin{pmatrix}
g_{11} & g_{12} & \gamma_{15}\\
g_{21} & g_{22} & \gamma_{25} \\
g_{31} & g_{32} & \gamma_{35} \\
g_{41} & g_{42} & \gamma_{45}\\
\gamma_{51} & \gamma_{52} & g_{55}
\end{pmatrix}
\quad \approx \quad
\begin{pmatrix}
I \\
W Z^{-1} \\
\rho_1 Z^{-1}
\end{pmatrix},
\begin{pmatrix}
I & 0 \\
V Y^{-1}& (\tau_2-WZ^{-1}\tau_1) a \\
0 & 1
\end{pmatrix},\label{bigcell}
\end{equation}
where 
\begin{eqnarray*}
&\rho_1=
\begin{pmatrix}
\gamma_{51} & \gamma_{52}
\end{pmatrix},\quad W=\begin{pmatrix}
g_{31} & g_{32} \\
g_{41} & g_{42}
\end{pmatrix},\quad Z=\begin{pmatrix}
g_{11} & g_{12} \\
g_{21} & g_{22}
\end{pmatrix}, \\&\tau_1=
\begin{pmatrix}
\gamma_{15} \\ \gamma_{25}
\end{pmatrix},\quad \tau_2=
\begin{pmatrix}
\gamma_{35} \\ \gamma_{45}
\end{pmatrix},\quad  d=(g_{55}-\nu Z^{-1}\mu_1)^{-1}\\&
V=W-g_{55}^{-1}\tau_2\rho_1,\quad Y=Z-g_{55}^{-1}\tau_1\rho_1.
\end{eqnarray*}
Finally the map $\pi$ in these coordinates  is given by:
$$
g \mapsto \left(W Z^{-1}, \rho_1 Z^{-1}, (\tau_2-WZ^{-1}\tau_1)
d\right).
$$
At this point one can compute the super Jacobian  and verify that
at the identity it is surjective.


Next, we are going to see how the big cell of the flag supermanifold
$\cF$ can be interpreted as the complex super
Minkowski space time, being the superflag its {\it superconformal
compactification}.

\medskip

The supergroup $G=\rSL_{4|1}$ is the complexification
of the real superconformal group.
The subgroup of  $G$
that leaves the big cell invariant is the set of matrices in $G$
of the form \begin{equation}\begin{pmatrix}L & 0& 0\\
NL&R&R\chi\\d\varphi
&0&d\end{pmatrix},\label{supergrouppoincare}\end{equation} with
$L, N, R$ being $2\times 2$ even matrices, $\chi$ and odd $1\times
2$ matrix, $\varphi$ a $2\times 1$ odd matrix and $d$ a scalar.
This is the complex Poincar\'{e} supergroup
and its action on
the big cell can be written as
\begin{eqnarray*}&&A\lra R(A+\chi\alpha)L^{-1}+N,\\&& \alpha\lra
d(\alpha +\varphi)L^{-1},\\&& \beta\lra
d^{-1}R(\beta+\chi).\end{eqnarray*} If the odd part is zero, then
the action reduces to the one of the classical
Poincar\'{e} group on the ordinary Minkowski
space (for more details see \cite{flv}).

\end{document}